\documentclass[12pt]{article}
\usepackage{amsmath}
\usepackage{amssymb}
\usepackage{url, enumerate}
\usepackage{mathtools}
\usepackage{framed}
\usepackage[framed]{ntheorem}
\newframedtheorem{frm-thm}{Theorem}

\newcommand{\eop}{\bigstar}  % end-of-proof

\newcommand{\otp}[1]{\hbox{otp($#1$)}}

\newenvironment{proof}{\noindent{\bf Proof.}}{\par\bigskip}

\newtheorem{THEOREM}{Theorem}[section]

\newtheorem{Conclusion}[THEOREM]{Conclusion}

\newtheorem{Hypothesis}[THEOREM]{Hypothesis}

\newtheorem{LEMMA}[THEOREM]{Lemma}

\newtheorem{Main Theorem}[THEOREM]{Main Theorem}
\newenvironment{main Theorem}{\begin{Main Theorem}} 
{\end{Main Theorem}}

\newtheorem{Theorem}[THEOREM]{Theorem}
\newenvironment{theorem}{\begin{Theorem}}{\end{Theorem}}

\newtheorem{Definition}[THEOREM]{Definition}

\newtheorem{Conventions}[THEOREM]{Conventions}

\newtheorem{Main Definition}[THEOREM]{Main Definition}
\newenvironment{main definition}{\begin{Main Definition}}
{\end{Main Definition}}

\newtheorem{Lemma}[THEOREM]{Lemma}
\newenvironment{lemma}{\begin{Lemma}}{\end{Lemma}}

\newtheorem{Notation}[THEOREM]{Notation}
\newenvironment{notation}{\begin{Notation}}{\end{Notation}}

\newtheorem{Convention}[THEOREM]{Convention}

\newtheorem{Note}[THEOREM]{Note}

\newtheorem{Observation}[THEOREM]{Observation}

\newtheorem{Remark}[THEOREM]{Remark}

\newtheorem{Question}[THEOREM]{Question}

\newtheorem{Main Fact}[THEOREM]{Main Fact}
\newenvironment{main Fact}{\begin{Main Fact}}{\end{Main Fact}}

\newtheorem{Fact}[THEOREM]{Fact}

\newtheorem{Subfact}[THEOREM]{Subfact}

\newtheorem{Claim}[THEOREM]{Claim}
\newenvironment{claim}{\begin{Claim}}{\end{Claim}}

\newtheorem{Main Claim}[THEOREM]{Main Claim}
\newenvironment{main claim}{\begin{Main Claim}}{\end{Main Claim}}

\newtheorem{Crucial Claim}[THEOREM]{Crucial Claim}
\newenvironment{crucial claim}{\begin{Crucial Claim}}{\end{Crucial Claim}}

\newtheorem{Subclaim}[THEOREM]{Subclaim}

\newtheorem{Sublemma}[THEOREM]{Sublemma}

\newtheorem{Corollary}[THEOREM]{Corollary}

\newtheorem{Example}[THEOREM]{Example}

\newtheorem{Problem}[THEOREM]{Problem}

\newtheorem{Proposition}[THEOREM]{Proposition}

\newtheorem{Conjecture}[THEOREM]{Conjecture}

\newtheorem{Discussion}[THEOREM]{Discussion}

\newenvironment{Proof of the Subfact}
{\noindent{\bf Proof of the Subfact.}}{\par\bigskip}

\newenvironment{Proof of the Theorem}
{\noindent{\bf Proof of the Theorem.}}{\par\bigskip}

\newenvironment{Proof of the Proposition}
{\noindent{\bf Proof of the Proposition.}}{\par\bigskip}

\newenvironment{Proof of the Conclusion}
{\noindent{\bf Proof of the Conclusion.}}{\par\bigskip}

\newenvironment{Proof of the Observation}
{\noindent{\bf Proof of the Observation.}}{\par\bigskip}

\newenvironment{Proof of the Fact}
{\noindent{\bf Proof of the Fact.}}{\par\bigskip}

\newenvironment{Proof of the Lemma}
{\noindent{\bf Proof of the Lemma.}}{\par\bigskip}

\newenvironment{Proof of the Claim}
{\noindent{\bf Proof of the Claim.}}{\par\bigskip}

\newenvironment{Proof of the Corollary}
{\noindent{\bf Proof of the Corollary.}}{\par\bigskip}

\newenvironment{Proof of the Subclaim}
{\noindent{\bf Proof of the Subclaim.}}{\par\medskip}

\newenvironment{Proof of the Main Claim}
{\noindent{\bf Proof of the Main Claim.}}{\par\bigskip}

\newenvironment{Proof of the Crucial Claim}
{\noindent{\bf Proof of the Crucial Claim.}}{\par\bigskip}

\newcount\skewfactor
\def\mathunderaccent#1#2 {\let\theaccent#1\skewfactor#2
\mathpalette\putaccentunder}
\def\putaccentunder#1#2{\oalign{$#1#2$\crcr\hidewidth
\vbox to.2ex{\hbox{$#1\skew\skewfactor\theaccent{}$}\vss}\hidewidth}}

\newcommand{\om}{\omega}

\author{Mirna D{\v z}amonja and Isa Vialard}
\title{On maximal order type of the lexicographic product}

\begin{document}
\maketitle
\begin{abstract} We give a self-contained proof of Isa Vialard's formula for $o(P\cdot Q)$ where $P$ and $Q$ are wpos. The proof introduces the notion of a cut of partial order, which might be of independent interest.
\end{abstract}

\section{Introduction} A well quasi order is a quasi order that has no infinite antichains or infinite decreasing sequences. For the purposes of this note, it is sufficient to work with well quasi orders that are actually partial orders, which we call well partial orders (wpo). A linearisation of a wpo is a linear order that extends the wpo but has the same underlying set. It follows that every linearisation of a wpo $P$ is a well order and has an ordinal order type. By a theorem of De Jongh and Parikh \cite{deJonghParikh} for every wpo $P$ there is the maximal ordinal which can be obtained as the order type of a linearisation of $P$, denoted $o(P)$ and called the maximal order type of $P$. In this note we are interested in the maximal order type of the orders obtained as the lexicographic product $P\cdot Q$ for some wpos $P$ and $Q$.

D{\v z}amonja, Schmitz and Schnoebelen incorrectly claimed in Lemma 4.4(1) of \cite{DzSS} that $o(P\cdot Q)=o(P)\cdot o(Q)$. Moreover, in \cite{DzSS} the statement was incorrectly attributed to Abraham and Bonnet. It turns out that while $o(P\cdot Q)=o(P)\cdot o(Q)$ is true if $o(Q)$ is a limit ordinal,
the statement is not true in general. The mistake was found by Harry Altman (private correspondance, March 2024), who provided a counter-example. Isa Vialard (March 2024, upcoming Ph.D. thesis) gave a correct formula for $o(P\cdot Q)$. This note presents his formula and a proof. Mirna D{\v z}amonja reformulated the original proof by Vialard into the more easily read terms used here, which we then jointly polished to obtain the final presentation.

The note contains Section \ref{intro} with the background, Section \ref{sec-cuts} with some new notions and Section \ref{Proof} with the actual proof.

\section{Background}\label{intro}

\subsection{Definition of $o(P)$} A partial order is {\em FAC} if it has no infinite antichains. It is {\em well founded} if it has no infinite decreasing sequence.
It is a {\em well partial order (wpo)} if it is both FAC and well founded.

\begin{claim}\label{maxelements} A FAC partial order has a finite (possibly 0) number of maximal elements.
\end{claim}

\begin{proof} Any two distinct maximal elements in a partial order are incompatible, hence the set of all maximal elements forms an antichain. Therefore, in
a FAC poset the set of all maximal elements is finite.
$\eop_{\ref{maxelements}}$
\end{proof}

\begin{claim}\label{linearisations} If $P$ is a wpo, then every one of its linearisations is a well order, so isomorphic to an ordinal.
\end{claim}

\begin{proof} Suppose for a contradiction that $L$ is a linearisation of $P$ but that $L$ has an infinite strictly decreasing sequence $\langle x_i:\,i<\omega\rangle$.
Then we note that for any $i<j$, either $x_i\bot_P x_j$ or $x_i>_P x_j$. Define a colouring of $[\omega]^2$ into the colours $\{0,1\}$ by
\begin{equation*}
c(i,j)=
\begin{cases} 0 & \text{ if } x_i\bot_P x_j \\
1 & \text{ if } x_i>_P x_j.
\end{cases}
\end{equation*}
Ramsey theorem gives an infinite mono-chromatic set for $c$. However, any 0-chromatic set for $c$ is an antichain and a 1-chromatic set a decreasing sequence, so none can be infinite in a wpo. Contradiction.
$\eop_{\ref{linearisations}}$
\end{proof}

For a wpo $P$ we define
\[
O(P)=\{\alpha: \,\mbox{there is a linearisation of }P\mbox{ of order type }\alpha\}.
\]

\begin{theorem}(de Jongh-Parikh \cite{deJonghParikh}) \label{deJonghParikh} For any wpo $P$, the set $O(P)$ has a maximal element.
\end{theorem}

\begin{notation}\label{o} If $P$ is a wpo, we denote $o(P)=\max O(P)$.
\end{notation}

\begin{claim}\label{mxelements} If $P$ is a wpo with $k$ maximal elements, then $o(P)=\delta + m$ for some limit ordinal $\delta$ (possibly 0)
and some natural number $m\ge k$. In addition, $k=0$ iff $m=0$.
\end{claim}

\begin{proof} We start with a couple of useful lemmas. For clarity in notation, in this proof we use $+_l$ to denote the lexicographic sum of two linear orders

\begin{lemma}\label{nomax} If $P$ is a wpo with no maximal elements, then every linearisation of $P$ has a limit order type.
\end{lemma}

\begin{proof} Suppose that $L$ is a linearisation of $P$ which has a successor order type, say $\beta+1$. In particular, $P$ is not empty. Without loss of generality, the elements of  $L$ and hence of $P$ are ordinals $\le \beta$, and hence $\beta$ is maximal in $L$. Since $P$ has no maximal elements, $\beta$ is not maximal in $P$ and hence there is $\alpha<\beta$ with $\beta<_P \alpha$. But then $\beta<_L \alpha$, in contradiction with $\beta$ being maximal in $L$.
$\eop_{\ref{nomax}}$
\end{proof}

\begin{lemma}\label{auxiliary} If $P$ is a wpo with at least one maximal element, $o(P)$ is not a limit ordinal.
\end{lemma}

\begin{proof} Suppose that $L$ is a linearisation of $P$ with order type $\delta$, a limit ordinal. Hence $\delta>0$ since $P\neq \emptyset$. Without loss of generality, the elements of $P$ and $L$ are the ordinals $<\delta$. Let $\beta<\delta$ be a maximal element of $P$. Then for every $\gamma\in (\beta,\delta)$ we have $\beta<_L \gamma$, which can only happen if $\beta\bot _P \gamma$.
It follows that $L\setminus \{\beta\} +_l \{\beta\}$
is a linearisation of $P$ of order type $\delta +1$, hence $\delta$ is not the maximal order type of a linearisation of $P$.
$\eop_{\ref{auxiliary}}$
\end{proof}

Now we prove the statement of the claim by induction on $k$. 
The case of $k=0$ follows directly from Lemma \ref{nomax}.

Suppose that we have proven the statement for $k \ge 0$ and let us prove it for $k+1$. Let $a$ be any maximal element of $P$ and let $Q=P\setminus \{a\}$.
Any maximal element of $P$ which is not equal to $a$ is necessarily a maximal element of $Q$, so $Q$ has at least $k$ maximal elements. 
By the induction hypothesis, $Q$ has a linearisation $L'$ of the type $o(Q')=\delta +m$ for some $m\ge k$.
Let $L=L' +_l\{a\}$, hence $L$ is a linearisation of $P$ of order type $\delta+(m+1)$ and therefore $o(P)\ge\delta+(m+1)$. We cannot have 
$o(P)\ge \delta+\omega$ since then we would also have $o(Q)\ge \delta+\omega$. Therefore, we have $o(P)=\delta+p$ for some $p\ge m+1\ge k+1$.
$\eop_{\ref{mxelements}}$
\end{proof}

\subsection{Products} Suppose that $P$ and $Q$ are two partial orders. Then $P\cdot Q$, {\em the lexicographic product} of $P$ along $Q$, also called {\em direct product} as in \cite{DzSS}, is
defined on the underlying set $P\times Q$ by letting
\[
(p_0,q_0)<(p_1,q_1) \mbox{ iff }q_0<_Q q_1 \vee (q_0=q_1 \wedge p_0<_P p_1).
\]

\begin{claim}\label{prod} If $P$ and $Q$ are FAC (wpo) posets, then so is $P\cdot Q$. 
\end{claim}

\begin{proof} Suppose that $P$ and $Q$ are FAC (wpo), but that $\langle (p_i, q_i):\,i<\omega\rangle$ is  an infinite antichain (decreasing sequence) in $P\cdot Q$. If
there is $q\in Q$ such that for an infinite set $A$ of indices $i$ in $\omega$ we have $q_i=q$, then the sequence
$\langle p_i: \, i\in A\rangle$ is an infinite antichain (decreasing sequence) in $P$, a contradiction. Hence there is an infinite set $B\subseteq \omega$ such that for any $i\neq j \in B$
we have $q_i\neq q_j$. 

In the case of FAC, since  $Q$ is FAC, there are $i<j$ such that either $q_i<_Q q_j$, so $(p_i,q_i)<_{P\cdot Q} (p_j,q_j)$, or 
$q_i>_Q q_j$, so $(p_i,q_i)>_{P\cdot Q} (p_j,q_j)$. In any case, $\langle (p_i, q_i):\,i<\omega\rangle$ is not an antichain. In the case of wpo, we have that there must be $i<j$ such that $q_i<_Q q_j$ and therefore $\langle (p_i, q_i):\,i<\omega\rangle$ cannot be strictly decreasing.
$\eop_{\ref{prod}}$
\end{proof}

Cantor Normal Formal for ordinals states that any ordinal can be represented uniquely as a sum of the form
\[
\omega^{\alpha_0}m_0+ \omega^{\alpha_1}m_1+\ldots + \omega^{\alpha_n}m_n,
\]
where $\alpha_0>\alpha_1>\ldots \alpha_n\ge 0$ are ordinals and $m_0, m_1,\ldots m_n$ are positive integers. The {\em Hessenberg sum} of two ordinals
$\alpha$ and $\beta$ is the ordinal $\alpha\oplus\beta$ whose Cantor normal form is the sum of Cantor normal forms of $\alpha$ and $\beta$, obtained by adding them as if they were polynomials.
We also define the {\em Hessenberg product} of two ordinals
$\alpha$ and $\beta$ as the ordinal $\alpha\otimes\beta$ whose Cantor normal form is the product of Cantor normal forms of $\alpha$ and $\beta$, obtained by multiplying them as if they were polynomials.

The following easily checked observation will be helpful in the main proof.

\begin{Observation}\label{helpfull} Suppose that $\alpha$ is an ordinal with the Cantor normal form
$\omega^{\alpha_0}m_0+ \omega^{\alpha_1}m_1+\ldots + \omega^{\alpha_n}m_n$ and $\beta$ is an ordinal with the Cantor normal form
$\omega^{\beta_0}o_0+ \omega^{\beta_1}o_1+\ldots + \omega^{\beta_l}o_l$. Further suppose that, $\alpha_n\ge \beta_0$.
Then 
\[
\alpha+\beta=\alpha\oplus\beta.
\]
\end{Observation}

\section{Cuts and isomorphic copies}\label{sec-cuts}
\begin{notation} If $P$ is a wpo and $L$ its linearisation of order type some ordinal $\alpha$, by renaming the elements of $P$ if necessary, we can assume that the elements of $P$ and $L$ as objects are ordinals in $\alpha$ and that the order on $L$ is the natural order on
$\alpha$. This implies that the order on $P$ satisfies
\[
\beta<_P \gamma \implies \beta <\gamma,
\]
since $L$ is a linearisation of $P$. We shall refer to this process of taking an isomorphic copy of $P$ as {\em without loss of generality, the elements of $P$ are ordinals in $\alpha$ and $P$ is a suborder of the ordinal order on $\alpha$}. Of course, neither $L$ nor $\alpha$ are unique, but once we fix $L$, they are, and this is the only type of the situation when we shall use this procedure.
\end{notation}

 Suppose that $P$ is a wpo and $P', P''\subseteq P$ are such that:
\begin{itemize}
\item $P'$ and $P''$ are disjoint,
\item the orders on $P'$ and $P''$ are inherited from $P$,
\item the underlying set of $P$ is the union of the underlying sets of $P'$ and $P''$ and
\item for all $p'\in P'$ and $p''\in P''$ either $p'\bot_P \,p''$ or $p'<_P p''$.
\end{itemize}
We say that $(P',P'')$ is a {\em cut} of $P$. 

\begin{claim}\label{cuts} Suppose that $(P',P'')$ is a cut of a wpo $P$. Then:
\begin{enumerate}
\item $P'$ and $P''$ are wpos.
\item $P$ is an augmentation of the disjoint union $P'\sqcup P''$ and is augmented by the lexicographic sum $P'+_l P''$.
\item $o(P')+o(P'')\le o(P)\le o(P')\oplus o(P'')$.
\item Suppose that there is a linearisation $L$ of $P$ of maximal order type such that 
\[
(\forall p'\in P')(\forall p''\in P'') p'<_L p''.
\]
Then $o(P)=o(P')+o(P'')$.
\end{enumerate}
\end{claim}

\begin{proof} (1) Follows as every infinite decreasing sequence (antichain) in $P'$ or $P''$ would give the rise to the same in $P$, since the order
on $P'$ and $P''$ is induced by that of $P$. 

{\noindent (2)} Follows since we cannot have $p'>_P p''$ for $p'\in P'$ and $p''\in P''$. 

{\noindent (3)} It is easily seen and proven in
Lemma 4.1 of \cite{DzSS} that $o(P'+_l P'')=o(P')+o(P'')$. Since every linearisation of $P'+_l P''$ is a linearisation of $P$,
we obtain $o(P')+o(P'')\le o(P)$. On the other hand, since $P$ is an augmentation of the disjoint union $P'\sqcup P''$, every linearisation of $P$
is also a linearisation of $P'\sqcup P''$ and hence $o(P)\le o(P'\sqcup P'')$. By Theorem 3.4 of \cite{deJonghParikh} we have 
$o(P'\sqcup P'')=o(P')\oplus o(P'')$, so proving the upper bound claimed.

{\noindent (4)} 
The equation is clearly true if either $P'$ or $P''$ is empty. Let us suppose that this is not the case.
Suppose that we have a linearisation $L$ of $P$ of maximal order type, say $\beta$ such that $(\forall p'\in P')(\forall p''\in P'') p'<_L p''$.
Without loss of generality, the elements of $P$ are the ordinals $<\beta$ and the order $<_P$ is a suborder of the ordinal order on $\beta$. 

We claim that this implies that  there is an ordinal $\alpha<\beta$ such that the elements of
$P'$ are the ordinals $<\alpha$ and the elements of $P''$ are those in $[\alpha, \beta)$. Let $\alpha=\min\{\gamma<\beta:\,\gamma\in P''\}$,
which is well defined since $P''\neq \emptyset$. Then clearly $[0,\alpha)\subseteq P'$. If $\gamma>\alpha$ and $\gamma\in P'$, then we have
$\alpha<_L\gamma$, in contradiction with $\alpha\in P''$ and $\gamma\in P'$.

Having established the above, we remark that $\alpha$ is a linearisation of $P'$, so $o(P')\ge \alpha$ and that $[\alpha, \beta)$ is a linearisation of 
$P''$, so $o(P'')\ge \otp{[\alpha, \beta)}$. Hence 
\[
o(P')+o(P'')\ge \alpha+  \otp{[\alpha, \beta)}=\beta=o(P).\eop_{\ref{cuts}}
\]
The other side of the inequality is already proved in (3).
\end{proof}

\section{$o(P\cdot Q)$}\label{Proof} In the light of \S\ref{intro}, one can legitimately ask the following question:

Suppose that $P$ and $Q$ are wpos, what is $o(P\cdot Q)$ in terms of $o(P)$ and $o(Q)$? As explained in the Introduction, an incorrect formula was claimed in Lemma 4.4(1) of \cite{DzSS} and the correct formula was given by Isa Vialard in March 2024. We now present that formula and its proof.

\begin{frm-thm}[Vialard]\label{Isasformula} Let $P$ and $Q$ be wpos and suppose that $Q$ has $k$ maximal elements. Hence $o(Q)=\delta +m$
for some limit ordinal $\delta$ and $m\ge k$.
Then
\[
o(P\cdot Q)=o(P)\cdot [\delta+(m-k)] +o(P)\otimes k.
\]
In particular, if $k=0$ (equivalently, $o(Q)$ is a limit)  then $o(P\cdot Q)=o(P)\cdot o(Q)$.
\end{frm-thm}

\begin{proof} If $P=\emptyset$ then $o(P)=0$, while $P\cdot Q=\emptyset$ and $o(P\cdot Q)=0$, which verifies the formula. Similarly, if $Q=\emptyset$
then $o(Q)=\delta=m=0$ and $P\cdot Q=\emptyset$, so the formula is verified.
Let us therefore assume that 
$P\neq \emptyset$ and $Q\neq \emptyset$.

Let us start by a general remark about the proposed formula.
Suppose that the Cantor normal form of $o(P)$ is $\omega^{\beta_0}k_0+\ldots + \omega^{\beta_l}k_l$ and denote 
$\sigma=\omega^{\beta_{1}}k_{1}+\ldots + \omega^{\beta_l}k_l$. Therefore we have $o(P)=\omega^{\beta_0}k_0+\sigma$, where $\sigma<
\omega^{\beta_0}$.
 Observe that for any $m\ge k$
\begin{equation*}
\begin{split}
o(P)\cdot(m-k)=\\
 (\omega^{\beta_0}k_0+\sigma )+ (\omega^{\beta_0}k_0+\sigma) +\ldots (\omega^{\beta_0}k_0+\sigma)=\\
\omega^{\beta_0}k_0 + (\sigma+ \omega^{\beta_0}k_0)+ (\sigma+ \omega^{\beta_0}k_0) +\ldots (\sigma+ \omega^{\beta_0}k_0) +\sigma= \\
\omega^{\beta_0}k_0+ \omega^{\beta_0}k_0 +\ldots  \omega^{\beta_0}k_0 +\sigma=\\
\omega^{\beta_0}k_0\cdot (m-k)+\sigma,
 \end{split}
 \end{equation*}
hence
\begin{equation*}
\begin{split}
o(P)\cdot [\delta+(m-k)] +o(P)\otimes k =\\
o(P)\cdot \delta + o(P)\cdot(m-k) + o(P)\otimes k =\\
o(P)\cdot\delta + \omega^{\beta_0}k_0\cdot (m-k)+\sigma
+  \omega^{\beta_0}k_0\cdot k+\sigma \otimes k =\\
o(P)\cdot\delta + \omega^{\beta_0}k_0\cdot (m-k)+ (\sigma
+  \omega^{\beta_0}k_0\cdot k)+\sigma \otimes k =\\
o(P)\cdot\delta + \omega^{\beta_0}k_0\cdot m
+\sigma \otimes k.
\end{split}
\end{equation*}

{\bf The induction.}
Let un now prove the formula. The proof is by induction on $o(Q)=\delta + m$.

{\bf The case of $Q$ finite, so $\delta=0$ and $m>0$}.
In this case $Q$ is a wpo of size $m$ with $1\le k\le m$ maximal elements
($k\ge 1$ since $Q$ is finite). 

Let $Q_\top$ be the restriction of $Q$ to the set of all maximal elements $q_1,\dots,q_k$ of $Q$, and let $Q_\bot=Q \setminus Q_\top$. 
Since $Q_\top$ is an antichain, we have that $P\cdot Q_\top$ is the disjoint union of
$\bigsqcup_{i\le k} P\times \{q_i\}$ and by Theorem 3.4 of \cite{deJonghParikh}, also cited as Theorem 4.2.(1) of \cite{DzSS}, we have that 
\begin{equation}\label{antichain}
o(P\cdot Q_\top)=\oplus_{i\le k} o(P)=o(P)\otimes k.
\end{equation}

We now obtain a lower bound on $o(P\cdot Q)$ using an induction on $m$.
Therefore
$o(P\cdot Q_\bot)\geq o(P)\cdot o(Q_\bot)$ by the induction hypothesis.
Furthermore, $(Q_\bot, Q_\top)$ is a cut of $Q$, thus $(P\cdot Q_\bot, P\cdot Q_\top)$ is a cut of $P\cdot Q$. 
Hence, according to Claim \ref{cuts}(3), and using that $o(Q_\bot)=m-k$, we obtain the lower bound

\begin{align*}
o(P\cdot Q)&\geq o(P\cdot Q_\bot)+ o(P\cdot Q_\top)\\
&\geq o(P)\cdot o(Q_\bot) + o(P)\otimes k\\
&= (\om^{\beta_0}k_0\cdot(m-k) + \sigma) +
(\om^{\beta_0}k_0\cdot k + \sigma\otimes k)\\
&= \om^{\beta_0}k_0\cdot(m-k) + (\sigma +
\om^{\beta_0}k_0\cdot k) + \sigma\otimes k\\
&=\om^{\beta_0}k_0\cdot m + \sigma\otimes k=o(P)\cdot [\delta+(m-k)] +o(P)\otimes k.\\
\end{align*}

 The proof of the upper bound is by induction on $k$, for all $m$ simultaneously.
 If $k=1$, then for any linearisation $L$ of $P\cdot Q$, 
 every element of $P\cdot Q_\bot$ is below every element of $P\cdot Q_\top$, thus applying Claim \ref{cuts}(4), we obtain
\begin{align*}
o(P\cdot Q)&= o(P\cdot Q_\bot)+o(P\cdot Q_\top)\\
&\leq o(P)\otimes (m-1) + o(P)\otimes 1 \mbox{ (by the induction hypothesis)}\\
&= (\om^{\beta_0}k_0\cdot(m-1) + \sigma\otimes (m-1))
+\om^{\beta_0}k_0 + \sigma\\
&= \om^{\beta_0}k_0\cdot m + \sigma,
\end{align*}
which matches the previous lower bound and so finishes the claim in this case.

If $k>1$ then let $Q''$ be the restriction of $Q$ to the subset containing $q_1$ and all the elements below $q_1$ that are not $\le_Q$ any other maximal element of $Q$, and let $Q'=Q\setminus Q''$. By definition, $Q'$ has $k-1$ maximal elements and $Q''$ has 1 maximal element. Let $m'=o(Q')$ and $m''=o(Q'')$, which by finiteness are in both cases the size of $Q'$ and $Q''$ respectively.
Then $m'\geq k-1$ and $m''\geq 1$, while $m=m'+m''$.
We know that $(Q', Q'')$ is a cut of $Q$, and therefore $(P\cdot Q', P\cdot Q'')$ is a cut of $P\cdot Q$. Hence applying Claim \ref{cuts}(3), 

\begin{align*}
o(P\cdot Q)&\leq o(P\cdot Q')\oplus o(P\cdot Q'')\\
&\overset{IH}{\le} [o(P)\cdot (m'-k+1) \oplus o(P)\otimes (k-1)] \oplus [o(P)\cdot (m''-1)\oplus o(P)]\\
&= (\om^{\beta_0}k_0\cdot m' + \sigma \otimes (k-1))
\oplus \om^{\beta_0}k_0\cdot m'' \oplus \sigma \\
&= \om^{\beta_0}k_0\cdot m + \sigma \otimes k,
\end{align*}
which matches the previous lower bound and so finishes the claim in this case.

{\bf The case of $o(Q)$ limit}. Suppose now that $m=0$, thus $Q$ has no maximal elements.
First observe that $o(P)$ and $o(Q)$ are augmentations of $P$ and $Q$, respectively, and therefore the linear order $o(P)\cdot o(Q)$ is a linearisation
of $P\cdot Q$. It follows that 
\begin{equation}\label{lowerbound}
o(P\cdot Q)\geq o(o(P)\cdot o(Q))= o(P)\cdot o(Q).
\end{equation}

For the upper bound, we fix a linearisation of maximal order type $L$ of $P\cdot Q$. By renaming we can assume that the elements of $P\cdot Q$ and $L$ are ordinals, so $L$ itself is an ordinal. For any fixed $\gamma<o(P\cdot Q)$, let us define $R_\gamma$ as the restriction of $P\cdot Q$ to $\gamma$. Without loss of generality, we may assume that $L$ maximises $o(R_\gamma)$, i.e., $o(R_\gamma)=\gamma$. Let $Q_\gamma$ be the smallest subset of $Q$ such that $R_\gamma\subseteq P \cdot Q_\gamma$. 
Note that $R_\gamma$ is downward-closed in $P\cdot Q$, and thus $Q_\gamma$ is downward-closed in $Q$. Since $o(Q)$ is limit, $Q$ has no maximal elements. Hence $Q_\gamma\neq Q$, because $R_\gamma$ is a strict downward-closed subset of $P\cdot Q$. Thus there is $x\in Q\setminus Q_\gamma$ such that $Q_\gamma\subseteq Q_{\not\geq x}$ since $Q_\gamma$ is downward-closed. Therefore $o(Q_\gamma)<o(Q)$.

% Note that $Q_\gamma\neq Q$, since $\gamma\notin P\cdot Q_\gamma$, but we shall moreover show that $o(Q_\gamma)<o(Q)$.

%Let $Q_{\gamma,\top}$ be the set of maximal elements of $Q_\gamma$ and $Q_{\gamma,\bot}=Q_\gamma\setminus Q_{\gamma,\top}$.
% Then we have
% $$o(P)\cdot o(Q_{\gamma,\bot})\leq o(P\cdot Q_{\gamma,\bot}) \leq o(R_\gamma)=\gamma,$$
%therefore $o(Q_{\gamma,\bot})<o(Q)$.  Since $Q_{\gamma,\top}$ is an antichain, it is finite. We have that $(Q_{\gamma,\bot},Q_{\gamma,\top})$ is a cut of $Q_\gamma$ so by Claim \ref{cuts}(3) it follows that  $o(Q_\gamma)\le o(Q_{\gamma,\bot})\oplus o(Q_{\gamma,\top})<o(Q)$ since $o(Q)$ is a limit.

Hence we can apply the induction hypothesis to obtain
$\gamma = o(R_\gamma)< o(P)\cdot o(Q)$ because the formula in Theorem \ref{Isasformula} is strictly increasing in $o(Q)$, independently on the number of maximal elements of $Q$. Since this is true for any $\gamma<o(P\cdot Q)$, we have proven that 
\begin{equation}\label{upperbound}
o(P\cdot Q)\leq o(P)\cdot o(Q).
\end{equation}
We obtain the desired equality by putting together the inequalities (\ref{lowerbound}) and (\ref{upperbound}).

{\bf The case of $o(Q)$ infinite and successor, so $\delta>0$ and $m>0$}. Let $o(Q)=\omega^{\alpha_0}m_0+\omega^{\alpha_1}m_1+\ldots + \omega^{\alpha_n}m_n$.
Hence $\alpha_0> 0$, $\alpha_n=0$ and $m_n=m>0$. Without loss of generality, $\delta+m$ is a maximal linearisation of $Q$ and elements of 
$Q$ are the ordinals in $\delta+m$.

Let $Q'$ be the restriction of $Q$ to $\delta$ and $Q''$ the restriction of $Q$ to 
$[\delta, \delta+m)$. 
Then:
\begin{enumerate}
\item $(Q',Q'')$ is a cut of $Q$ satisfying the requirement of Claim \ref{cuts}(4), hence $o(Q)=o(Q')+o(Q'')$.
\item As $Q''$ is finite of size $m$, clearly  $o(Q'')=m$.
\item\label{treci} $\delta$ is a linearisation of $Q'$ of order type $\delta$ and replacing $\delta$ by any longer linearisation of $Q'$, say $\sigma >\delta$,
would give a linearisation of $Q$ of order type $\ge \sigma +m\ge (\delta+1)+m>\delta + m$, a contradiction. Hence $o(Q')=\delta$.
\item\label{cet} $Q'$ has no maximal element, as it has a linearisation $L'$ with no
maximal element (see Lemma \ref{nomax}). 
\item $Q''$ has $k$ maximal elements. This is so since any maximal element of $Q$ which is in $Q''$ is also maximal in $Q''$, but by (\ref{cet}) above, all maximal elements of
$Q$ are in $Q''$. Hence $Q''$ has at least $k$ maximal elements. But if it had a maximal element $q$ which is not a maximal element in $Q$, then there
would be $p\in Q'$ with $p>_P q$, which is a contradiction with $(Q',Q'')$ being a cut.
\end{enumerate}
Now notice that $(P\cdot Q', P\cdot Q'')$ is a cut of $P\cdot Q$, since if $(p',q')\in P\cdot Q'$ and $(p'',q'')\in P\cdot Q''$ were to be such that 
$(p',q')\ge_{P\cdot Q} (p'',q'')$, then, since $Q'\cap Q''=\emptyset$, we would have that $q'>_Q q''$, which is a contradiction with $(Q',Q'')$ being a cut of
$Q$.
Then: 
\begin{description}
\item{(a)} $(P\cdot Q', P\cdot Q'')$ is a cut of $P\cdot Q$ so by Claim \ref{cuts}(3) we have
\[
o(P\cdot Q')+ o(P\cdot Q'')\le o(P\cdot Q)\le o(P\cdot Q')\oplus o(P\cdot Q'').
\]
\item {(b)} $o(P\cdot Q')=o(P)\cdot \delta$  by the induction hypothesis, since $o(Q')$ is a limit $\delta<\delta+m$.
\item {(c)} $o(P\cdot Q'')=o(P)\cdot(m-k) + o(P)\otimes k$, by the induction hypothesis that we proved in the case of $Q$ finite.
\item{(d)} $o(P\cdot Q')+ o(P\cdot Q'')= o(P)\cdot \delta + o(P)\cdot (m-k) +o(P)\otimes k=o(P) \cdot [\delta+(m-k)])+o(P)\otimes k$.
\item{(e)} $o(P\cdot Q')\oplus o(P\cdot Q'')= o(P)\cdot \delta \oplus  [o(P)\cdot (m-k) +o(P)\otimes k]$. 

Let the Cantor normal form of $o(P)$ be $\omega^{\beta_0}o_0+ \omega^{\beta_1}o_1+\ldots + \omega^{\beta_l}o_l$.
By inspection, we have that the
smallest exponent in the Cantor normal form of $o(P)\cdot \delta$ is at least $\beta_0$ and that the greatest 
exponent in the Cantor normal form of $o(P)\cdot (m-k) +o(P)\otimes k$ is $\beta_0$. Hence, By Observation \ref{helpfull}, we have
$o(P)\cdot \delta \oplus  [o(P)\cdot (m-k) +o(P)\otimes k]=o(P)\cdot \delta \oplus  [o(P)\cdot (m-k) +o(P)\otimes k]$. 
\end{description}
By Claim \ref{cuts}(3) we have 
\[
o(P\cdot Q)= o(P)\cdot [\delta+(m-k)] +o(P)\otimes k,
\]
as required.
$\eop_{\ref{Isasformula}}$
\end{proof}

\bibliographystyle{plain}
\bibliography{../bibliomaster}

\end{document}